\documentclass[10pt,twocolumn,twoside]{IEEEtran}
\usepackage{multirow}
\usepackage{xspace}

\newcommand{\ie}{i.e.\xspace}
\usepackage{mathtools}
\usepackage[dvipsnames]{xcolor}

\DeclarePairedDelimiter\abs{\lvert}{\rvert}%
\newcommand{\calE}{\mathcal{E}}

\newcommand{\calG}{\mathcal{G}}

\newcommand{\calP}{\mathcal{P}}

\newcommand{\calV}{\mathcal{V}}

\newcommand{\bbR}{\mathbb{R}}

\newcommand{\defeq}{\vcentcolon=}

\usepackage{booktabs}
\usepackage{graphicx}
\usepackage{mathtools}
\usepackage{booktabs}
\usepackage{graphicx}
\usepackage{mathtools}

\newcommand{\binaryset}{\{0, 1\}}

\newcommand{\GcalE}{G^{\calE}}

\newcommand\clearrow{\global\let\rowmac\relax}
\clearrow

\usepackage{array}
\newcolumntype{H}{>{\setbox0=\hbox\bgroup}c<{\egroup}@{}}

\newcommand{\alen}[1]{}

\newcommand{\tppgamma}{TP-NC}
\newcommand{\tppfd}{TP-PFD}

\usepackage{optidef}
\usepackage{tabularx}
\usepackage{multirow}
\usepackage{graphicx}
\usepackage{booktabs}
\usepackage{bm}
\usepackage{hyperref}
\usepackage[capitalise]{cleveref}
\crefname{equation}{}{}
\usepackage{amsmath, amsthm}
\usepackage{amssymb}
\usepackage{siunitx}
\usepackage{placeins}
\usepackage[caption=false,font=footnotesize]{subfig}
\graphicspath{{figs/}}

\newcommand{\edit}[1]{#1}
\ifCLASSINFOpdf
\else
\fi
\begin{document}
\title{Mixed-integer linear programming approaches \\for tree partitioning of power networks}

\author{Leon~Lan~and~Alessandro~Zocca% <-this % stops a space
% \thanks{?}% <-this % stops a space
% \thanks{Manuscript received April 19, 2005; revised August 26, 2015.}
}


\markboth{Lan and Zocca: MILP approaches for tree partitioning of power networks}{}



\maketitle
\begin{abstract}
  In transmission networks, power flows and network topology are deeply intertwined due to power flow physics. Recent literature shows that a specific more hierarchical network structure can effectively inhibit the propagation of line failures across the entire system.
  In particular, a novel approach named \textit{tree partitioning} has been proposed, which seeks to bolster the robustness of power networks through strategic alterations in network topology, accomplished via targeted line switching actions.
  Several tree partitioning problem formulations have been proposed by considering different objectives, among which power flow disruption and network congestion. Furthermore, various heuristic methods based on a two-stage and recursive approach have been proposed. 
  The present work provides a general framework for tree partitioning problems based on mixed-integer linear programming (MILP). In particular, we present a novel MILP formulation to optimally solve tree partitioning problems and also propose a two-stage heuristic based on MILP. We perform extensive numerical experiments to solve two tree partitioning problem variants, demonstrating the excellent performance of our solution methods. Lastly, through exhaustive cascading failure simulations, we compare the effectiveness of various tree partitioning strategies and show that, on average, they can achieve a substantial reduction in lost load compared to the original topologies.
\end{abstract}

\begin{IEEEkeywords}
Tree partitioning,  controlled islanding, mixed-integer linear programming, cascading failures, line switching
\end{IEEEkeywords}

%%% Local Variables: ***
%%% mode:latex ***
%%% TeX-master: "main.tex"  ***
%%% End: ***

\section{Introduction}
\edit{Power systems are the backbone of modern society, whose functioning nowadays significantly relies on a continuous and stable power supply.} Consequently, their reliability is of paramount importance: any disruption or instability in the power supply can have severe economic and societal implications.
Historical power blackouts have shown that transmission line failures play a critical and often initiating role in the unfolding of cascading failures \cite{Vaiman2012,Bienstock2015}. \edit{The complex interplay of grid topology and power flow physics gives rise to intricate failure patterns, often exhibiting a non-local propagation of line failures~\cite{Bernstein2014,Dobson2016,Hines2017,Witthaut2015,Guo2020}.}

% Emergency measures and zooming in on CI 
Several emergency measures can be deployed to support frequency and voltage stability and mitigate cascading failures in transmission power systems, among which dynamic line rating, generation re-dispatch, load shedding, and the activation of fast-reacting reserves. \edit{When these measures fail, a last-resort emergency measure to stop cascading failures is \textit{controlled islanding}, which splits the network into separate connected components called ``islands''~\cite{Senroy2006}. By confining failures within isolated islands, this approach helps first to stabilize and later to restore the system more effectively.} 

\begin{figure}[!bt]
    \centering
    % Trim: left bottom right top
    \includegraphics[width=\linewidth%,trim={1.8cm 1.4cm 1.9cm 2.5cm}
    ]{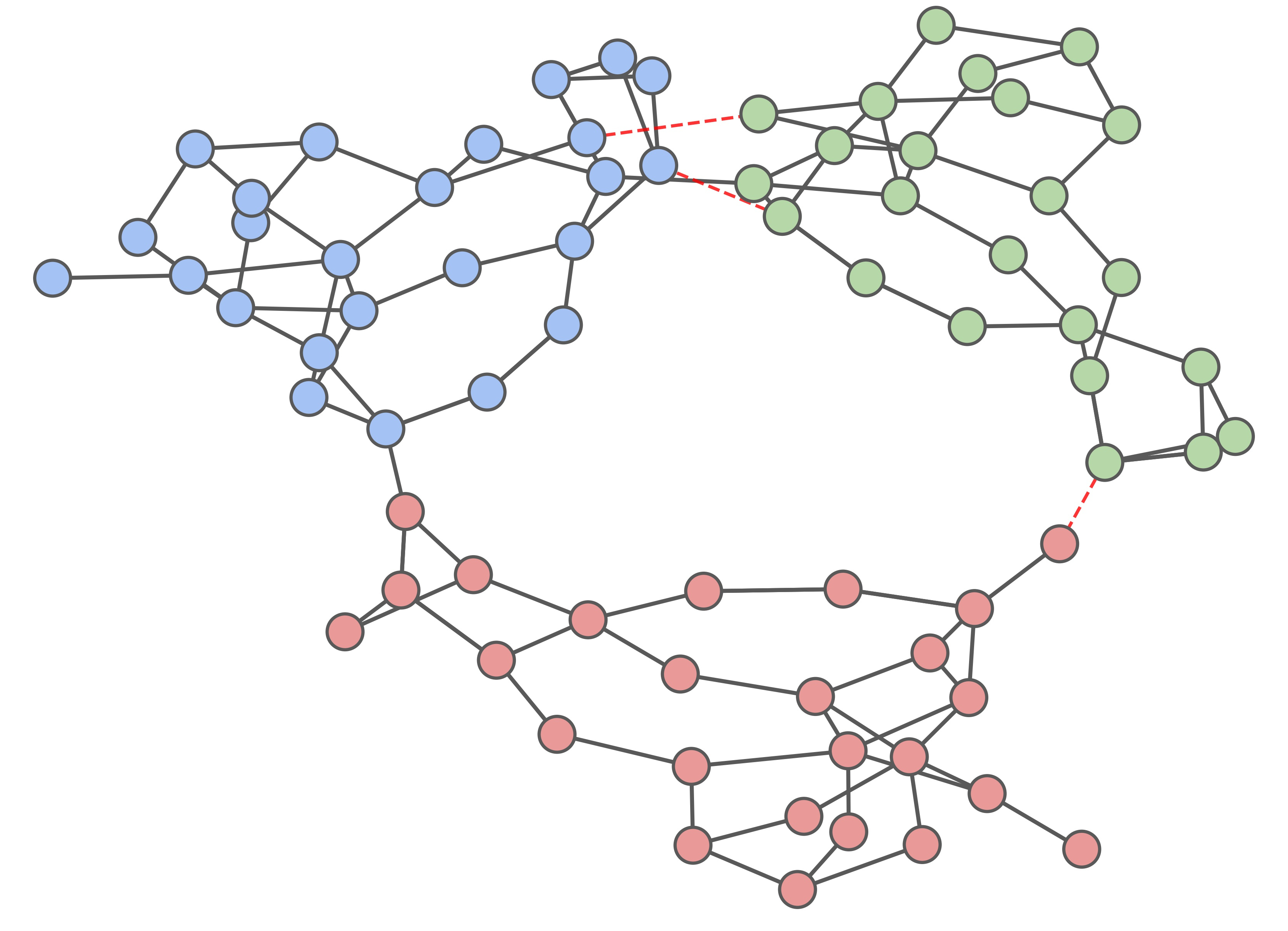}
    \caption{
    Tree partitioning of the IEEE-73 network with provable line failure localization properties. 
    The three clusters connected in a tree-like manner have been obtained by switching off only three transmission lines (in red).
    }
    \label{fig:tp_intro}
\end{figure}

However, despite the widespread attention that controlled islanding has received in the literature \cite{QianchuanZhao2003,Ding2013,Trodden2013,Sanchez-Garcia2014,Quiros-Tortos2015}, its practical implementation has been limited \cite{BialekJanusz2021}.
% TP as an alternative to CI
A recent paper \cite{BialekJanusz2021} proposes an alternative less drastic emergency measure, named \textit{tree partitioning}, that takes advantage of the same network flexibility in terms of line switching actions. 
\edit{Line switching is a consolidated electricity grid management paradigm that received a lot of attention also in the optimal transmission switching literature~\cite{Salkuti2018,Hedman2011}, where line switching actions aim at maximizing economic efficiency of generation dispatch.}

Unlike controlled islanding, which creates disconnected components, tree partitioning modifies the power network to obtain a hierarchical structure with clusters interconnected in a tree-like topology (see Figure~\ref{fig:tp_intro}).
As shown in~\cite{BialekJanusz2021,Zocca}, any line failure in a network after tree partitioning causes a power flow redistribution (and hence possible subsequent failures) only within the affected cluster.
\edit{Thus, a tree-partitioned network has the same ability to locally confine line failures and mitigate cascading failures as the corresponding islanding strategy, while offering several advantages:} (i) avoiding load shedding, (ii) reduced impact on the surviving network (needing fewer switching actions and less power flow redistribution), (iii) no need for generator adjustments and (iv) no need for re-synchronization when resuming normal operations \cite{BialekJanusz2021}.

Different problem formulations have been proposed to determine the best tree partitioning strategy, focusing on minimizing power flow disruption \cite{BialekJanusz2021} or network congestion \cite{Zocca}.
\edit{The first formulation aims to} minimize the impact of switching actions on the transient stability of the network, considered also often in the controlled islanding literature \cite{Ding2013,Quiros-Tortos2015,Demetriou2016}. 
\edit{In the second formulation, the goal is to minimize network congestion, ensuring minimal line overloading after partitioning.
These problem variants have been addressed with two-stage approaches, with the network congestion variant also explored via recursive algorithms. 
However, existing methods are heuristic and provide no guarantees for optimal switching actions in tree partitioning.}

This paper introduces a general framework for tree partitioning problems and describes approaches based on mixed-integer linear programming (MILP) to solve them.
More specifically, the main contributions of this paper are as follows:
\begin{itemize}
    \item We develop a novel MILP formulation that unifies existing problem variants, obtaining for the first time optimal tree partitioning solutions.
    \item We propose a heuristic two-stage approach based on MILP to improve the scalability for large networks. 
    \item We extensively compare the performance of our proposed solution methods on a large collection of test cases, demonstrating that the exact approach achieves significantly better solutions, while the heuristic method obtains good solutions with considerably lower runtimes.
    \item We compare the performance and effectiveness of various tree partitioning strategies through extensive cascading failure simulations, showing that, on average, tree partitioning can achieve a substantial reduction in lost load compared to the original topologies.
\end{itemize}
 
The rest of the paper is structured as follows. Section~\ref{sec:preliminaries} presents definitions and mathematical preliminaries. Section~\ref{sec:problem} introduces a general tree partitioning formulation and its problem variants. In Section~\ref{sec:methods}, we describe solution methods to solve tree partitioning problems. Section~\ref{sec:experiments} presents numerical experiments of the tree partitioning problems, as well as cascading failure simulations. Finally, Section~\ref{sec:conclusion} concludes the paper.

%%% Local Variables: ***
%%% mode:latex ***
%%% TeX-master: "main.tex"  ***
%%% End: ***
\section{Preliminaries}
\label{sec:preliminaries}
The goal of this section is to describe the power network model we consider in this work (cf.~\cref{sec:power-network-model}), introduce the notion of tree partitions (cf.~\cref{sec:tree-partitions}), and review the failure localization properties that power networks have after tree partitioning (cf.~\cref{sec:failure-localization}).

\subsection{Power network model}
\label{sec:power-network-model}
We model an electrical transmission network as a connected, directed graph $G=(V,E)$, where $V=\{1, 2, \dots, n\}$ is the set of vertices (buses) and $E \subseteq \{(i, j): i, j \in V\}$ is the set of edges (transmission lines). Let $n$ and $m$ denote the number of buses and lines, respectively. Each bus \(i \in V\) has a \textit{net power injection} \(p_i\), where \(p_i > 0\) is interpreted as \emph{injected} power and \(p_i \leq 0\) as \emph{consumed} power. Each line $(i,j) \in E$ has a capacity \(c_{ij} > 0\), denoting its rating, i.e., the maximum power that the line can safely carry.

In this paper, we consider a lossless DC power flow model in which generation always matches demand, i.e., $\sum_{i=1}^{n} p_i = 0$. We refer to any such vector $\bm{p}$ of power injections as \textit{balanced}. Let $\theta_{i} \in [\theta^{\min}, \theta^{\max}]$ denote the \textit{phase angle} of bus $i$, where $\theta^{\min}, \theta^{\max}  \in \mathbb{R}$ denote the minimum and maximum phase angles, respectively.
For each line $(i,j)$, let $f_{ij} \in \bbR$ denote the active power flow and let \(b_{ij} > 0\) denote the \textit{line susceptance}. Given a vector of power injections $\bm{p}$, the corresponding line flows $\bm{f}$ and phase angles $\bm{\theta}$ are obtained by solving the \textit{DC power flow equations}:
\begin{subequations}
\label{eq:DCPF}
\begin{align}
p_i &= \sum_{j: (i, j) \in E} f_{ij} - \sum_{j: (j, i) \in E} f_{ji}, & \forall \, i \in V, \label{eq:DCPF-1}\\
  f_{ij} &=  b_{ij}(\theta_i - \theta_j), & \forall \, (i, j) \in E. \label{eq:DCPF-2}
\end{align}
\end{subequations}
Equation~\eqref{eq:DCPF-1} ensures flow conservation, and \eqref{eq:DCPF-2} captures the flow dependency on susceptances and angle differences. The DC power flow equations \eqref{eq:DCPF} admit a unique power flow solution $\bm{f}$ for each balanced injection vector \(\bm{p}\).

\subsection{Tree partitions}
\label{sec:tree-partitions}
\edit{To introduce the notion of tree partition, we first quickly revise some key graph theory concepts. Given a graph $G=(V,E)$, an edge $e \in E$ is called a \textit{cut-edge} or \textit{bridge} if its removal would disconnect the graph.} A \textit{$k$-partition} of $G$ is the collection $\calP = \{ \calV_1, \calV_2, \dots, \calV_k \}$ of $k$ non-empty, disjoint vertex sets $\calV_1, \calV_2, \dots, \calV_k$, called \textit{clusters}, such that $\bigcup_{i=l}^k \calV_l = V$.
We denote by $K = \{1, 2, \dots, k\}$ the set of cluster indices.
Given a partition $\calP$, an edge $(i,j) \in E$ is called an \textit{internal edge} if both $i$ and $j$ belong to the same cluster and \textit{cross edge} otherwise. The set of cross edges is denoted by $E_{C}(\calP)$. 

\begin{figure}[!h]
    % [trim = {left, bottom, right, top}, clip]
    \centering
    \subfloat[]{{\includegraphics
    [width=0.48\linewidth,trim={5cm 1cm 5cm 1cm},clip]
    {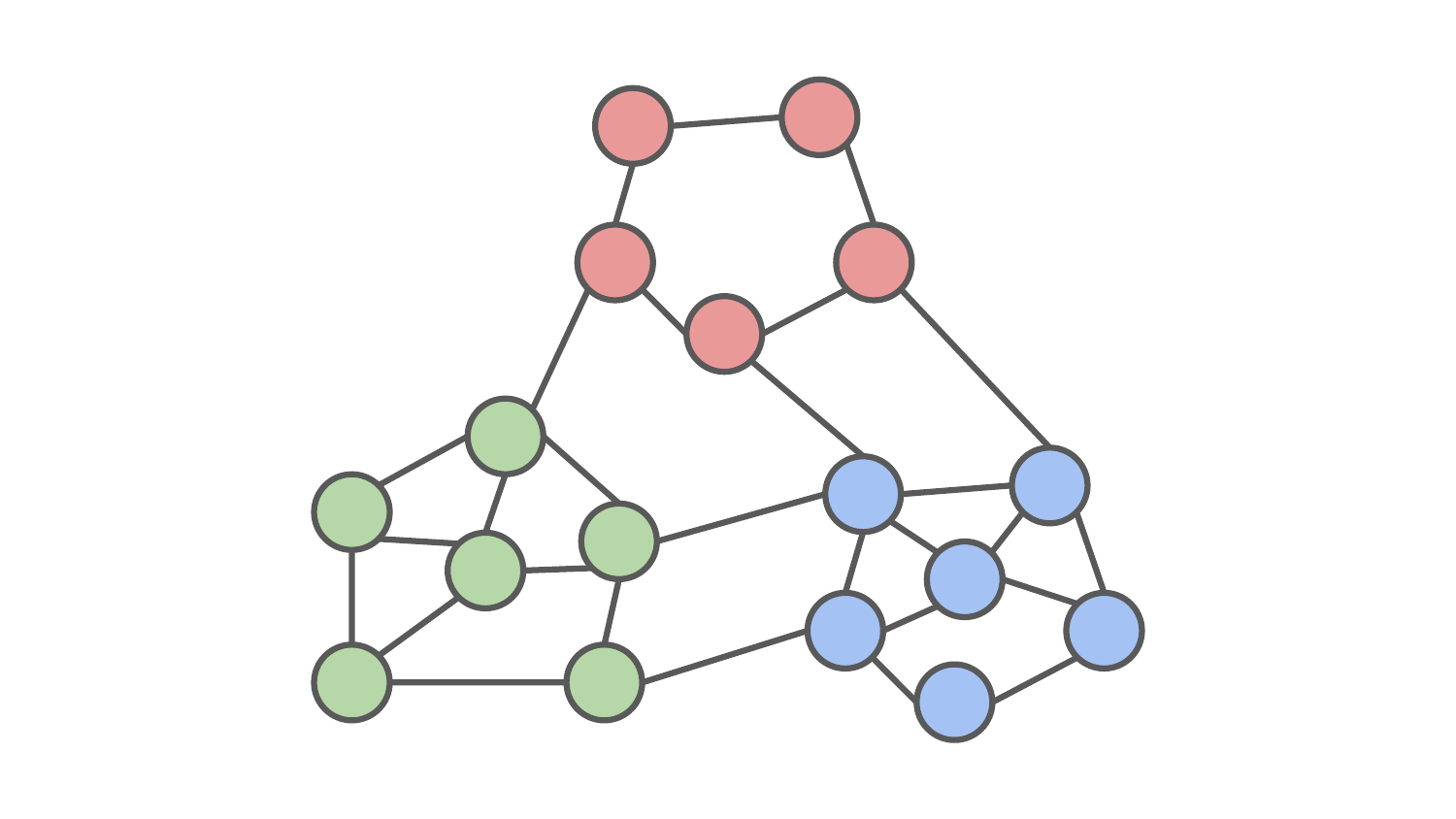} }}%
    \hfill
    \subfloat[]{{\includegraphics
    [width=0.48\linewidth,trim={5cm 1cm 5cm 1cm},clip]
    {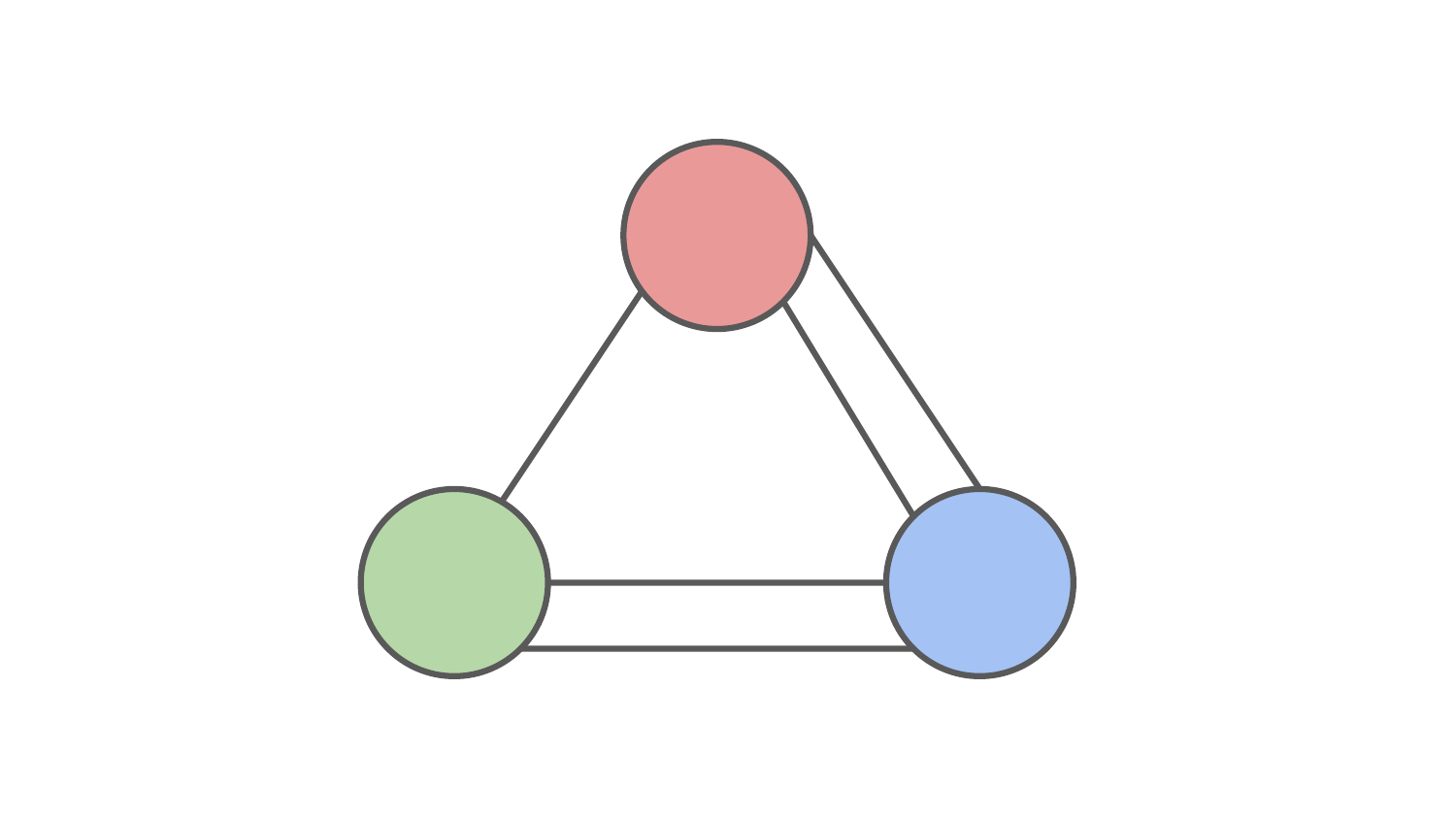} }}%
    \caption{(a) A graph $G$ with colors representing a partition $\calP$ into three clusters. (b) The reduced graph $G_\calP$ that corresponds to the partition in (a).}%
    \label{fig:reduced_graph}
\end{figure}

Given a partition $\calP$ of the graph $G$, the corresponding \textit{reduced graph} $G_\calP=(K, E_{C}(\calP))$  is the graph whose vertices are the clusters in $\calP$ and where an edge is drawn for each cross edge connecting two different clusters (see Figure~\ref{fig:reduced_graph}). Note that it is possible for the reduced graph to have multiple edges between two vertices, and thus to be a multigraph. 
We say that $\calP$ is a \textit{tree partition} of $G$ if the corresponding reduced graph $G_\calP$ is a tree. \edit{The cross edges of a tree partition are automatically cut-edges for the graph $G$.} %If $\calP$ is a tree partition, then every cross edge is a \textit{bridge} for the graph $G$. 
% \leon{Cut edges are not defined. Do we assume this is common graph theory terminology?}

In general, a graph $G$ admits multiple tree partitions, including the trivial one where all nodes are in the same cluster. %as shown in Figure~\ref{fig:tree-partition}~and~\ref{fig:bbd}. 
\edit{There is a unique tree partition that is maximal by inclusion~\cite{Guo2018}. We refer to this as the \textit{irreducible tree partition}.} 
%or as the \textit{bridge-block decomposition} of the graph. This name is due to the fact that the cross edges corresponding to the irreducible tree partition are precisely all the \edit{\textit{cut-edges}} of the original graph.}

\subsection{Failure localization}
\label{sec:failure-localization}

When a transmission line outage occurs (either due to failure or intentional disconnection), the power flow it carried is redistributed \textit{globally} across the remaining lines according to power flow physics. \edit{This redistribution can cause line flows to increase, decrease, or even reverse direction.} Some of the remaining lines can then become overloaded and are deactivated by automatic safety relays, leading to \textit{failure propagation} in the network. This phenomenon is not localized: as shown in~\cite{Kinney2005} using real-world data, line failures can trigger subsequent line failures even very far away from the initial contingency.

It was shown in \cite{BialekJanusz2021, Zocca, Guo2018, Guo2020} that line failures do not propagate across cut-edges, \edit{meaning that, if one considers the clusters determined by the irreducible tree partition, line failures only impact lines flows within the same cluster.}
These results have been made rigorous in~\cite{Guo2020, Guo2020a} under the DC power flow approximation using the so-called \textit{line outage distribution factors (LODF)}~\cite{Guler2007a, Soltan2017a}, \edit{a standard tool in the power systems literature to compute the post-contingency line flows. 
In particular, \cite{Guo2020, Guo2020a} show that finer tree partitions of the same network lead to sparser LODF matrices. Although a formal proof of the same result for the AC power flow model is lacking,} there is substantial numerical evidence from AC simulations~\cite{Guo2020a} that tree partitions effectively localize an overwhelming majority of line failures also in this setting.

\edit{Most power networks, however, have a highly meshed structure, making their irreducible tree partition trivial and thus very prone to non-local line failure propagation \cite{Zocca}. In the next section, we formulate an optimization problem that aims to enhance the robustness of a power network against cascading line failures by refining its irreducible tree partition using line-switching actions.}

%%% Local Variables: ***
%%% mode:latex ***
%%% TeX-master: "main.tex"  ***
%%% End: ***
\section{Tree partitioning problem formulations}
\label{sec:problem}

In light of the strong interplay between network topology and line failure propagation described in~\Cref{sec:failure-localization}, \edit{tree partitioning strategies~\cite{BialekJanusz2021, Zocca, Guo2018} aim to slightly modify the topology of a given power network to enhance failure localization using line-switching actions.

Given a power network $G=(V, E)$, we want to identify the optimal subset $\calE \subset E$ of transmission lines to switch off to maximize the failure localization properties of the post-switching network $\GcalE=(V, E \setminus \mathcal{E})$ in terms of its irreducible tree partition. 
This involves identifying a partition $\calP$ and a set of lines $\calE$ such that, once the lines in $\calE$ are deactivated using switching actions, $\calP$ becomes a tree partition for the post-switching network $\GcalE$. We will refer to $\GcalE$ as the (resulting) \textit{tree-partitioned network}.}

Various formulations of the tree partitioning problem have been proposed, but the core structure is the same and can be described as follows. Given a power network $G=(V,E)$ and $k \geq 2$ coherent generator groups $\mathcal{G}_1, \mathcal{G}_2, \dots, \mathcal{G}_k \subset V$, the goal is to identify a $k$-partition $\calP = \{ \mathcal{V}_1, \mathcal{V}_2, \dots, \mathcal{V}_k \}$ and a subset of lines $\calE \subset E_C(\calP)$ that minimize a specific risk function $r(\calE)$ subject to the following constraints:
\begin{itemize}
    \item[(a)] the post-switching network $\GcalE=(V, E \setminus \mathcal{E})$ is connected;
    \item[(b)] $\calP$ is a tree partition of $G^\calE$;
    \item[(c)] the partition is such that coherent generators belong to the same cluster, i.e., $ \mathcal{G}_l \subseteq \mathcal{V}_l,$ for every $l \in K := \{1,\dots,k\}$.
\end{itemize}

Constraint (c), which we refer to as the \textit{coherent generator grouping constraint}, arises naturally when thinking of tree partitioning as an emergency measure against cascading failures, as proposed by~\cite{BialekJanusz2021}. In this context, tree partitioning inherits several design principles behind islanding schemes~\cite{Ding2013, Quiros-Tortos2015, Demetriou2016}.
In particular, it is key to minimize the impact on the \textit{transient stability}, that is, the ability of the power network to maintain frequency synchronization when subjected to a severe disturbance \cite{Machowski2020}. 
It is therefore important to maintain \textit{generator coherency}: synchronous generators 
should preferably be grouped together in the same cluster. In fact, non-coherent generator groups might lead to generator tripping and, eventually, the collapse of the network \cite{Ding2013}.
% The coherent generator grouping constraint also has the advantage of reducing the complexity of the problem since otherwise the number of possible pairs $(\calP,\calE)$ would grow exponentially in the network size.
Various other additional constraints can be added to this core problem, capturing either desired features of the partition $\calP$, or physical properties of the power network.

Several choices can also be made for the objective function $r(\calE)$, but they all try to capture the risk of removing from service the selected lines in $\calE$ by quantifying the impact of such switching actions. \edit{Deactivating lines redistributes power flows (as discussed in~\cref{sec:failure-localization}) and may influence system stability. Thus, the optimal tree partition and switching actions depend on the specific application.} In the next two sections, \Cref{sec:transient-stability,sec:network-congestion}, we explore two variants of the tree partitioning problem introduced by \cite{BialekJanusz2021,Zocca}.

\subsection{Minimizing power flow disruption}
\label{sec:transient-stability}
In controlled islanding, the impact of switching actions on the rest of the network is usually quantified using the \textit{power flow disruption}, defined as the total sum of absolute flows on the switched-off lines. For a tree-partitioned network, constraint (a) ensures that the network remains connected and thus there is no actual power disruption -- one of its key advantages over islanding. Nevertheless, it is still reasonable to consider the same quantity, i.e., the sum of the absolute flows on the lines in $\calE$ selected to be switched off, as a proxy for the impact of line-switching actions on the network, namely
\[
     r(\calE) \defeq \sum_{(i, j) \in \calE} \abs{f_{ij}}.
\]

Using this objective function, we can now formulate precisely the \textit{tree partitioning problem considering power flow disruption} (TP-PFD). Given a power network $G=(V,E)$ and $k$ coherent generator groups $\calG_1, \calG_2, \dots, \calG_k \subset V$, the goal of TP-PFD is to identify a $k$-partition $\calP = \{\calV_1, \calV_2, \dots, \calV_k\}$ and a set of lines $\mathcal{E} \subset E$ to switch off such that the power flow disruption $\sum_{(i, j) \in \mathcal{E}}|f_{ij}|$ is minimized while constraints (a), (b), and (c) are satisfied:
\begin{align*}
    \min \; & \sum_{(i, j) \in E} |f_{ij}|,\\
    \text{s.t. } & \text{(a), (b), (c)}.
\end{align*}

\subsection{Minimizing network congestion}
\label{sec:network-congestion}
Tree partitioning has been proposed not only as an emergency measure but also as a tool in the planning phase to obtain a network that is simultaneously more reliable and requires less effort to mitigate disturbances~\cite{Zocca, Liang2022a}. In this context, the impact of line-switching actions can be quantified differently. 
\edit{Specifically, \cite{Zocca} proposes to look at the congestion on the lines in the post-switching network $\GcalE$, after the power flow redistribution induced by the selected switching actions}. %as a result of flow redistribution after the switching actions.

We define the \textit{congestion level} of a single line $(i,j) \in E$ as the non-negative ratio $\abs{f_{ij}} \mathbin{/} c_{ij}$, and we say that a line $(i, j) \in E$ is congested if its congestion level exceeds one. For a given power network $G=(V,E)$, we define the \emph{network congestion} as
\[
  \Gamma(G) \defeq \max_{(i, j) \in E} \abs{f_{ij}} \mathbin{/} c_{ij}.
\]
i.e., as the maximum congestion level \edit{over all network lines. If $\Gamma(G) \le 1$, then the network has no congested lines.}

The authors in~\cite{Zocca} propose a variant of the tree-partitioning problem in which the goal is to minimize the impact of switching actions on the network congestion level. In our framework, this amounts to considering as the objective function the post-switching network congestion $\Gamma(\GcalE)$, i.e.,
\[
    r(\calE) \defeq \Gamma(\GcalE).
\]

Here $\Gamma(\GcalE)$ is calculated on the network $G^\calE$ assuming that the power injections $\bm{p}$ are unchanged, and the post-switching power flows are recalculated using DC power flow equations \eqref{eq:DCPF}.
\edit{The rationale behind this problem formulation} is that, when tree partitioning a network, it is desirable to switch off a set of lines $\calE$ so that the network congestion on the post-switching network $\GcalE$ is as low as possible and, in particular, remains below one. 
\edit{However, due to the complex interplay between power network topology and power flow physics, it is hard to ensure no congestion upfront.}
Nevertheless, slight network congestion may still occur, but it can often be alleviated with remedial actions~\cite{BialekJanusz2021}.

We can now formulate the \textit{tree partitioning problem considering network congestion} (\tppgamma). Given a power network $G=(V,E)$ and $k$ coherent generator groups $\calG_1, \calG_2, \dots, \calG_k \subset V$, the goal of \tppgamma~is to identify a $k$-partition $\calP = \{\calV_1, \calV_2, \dots, \calV_k\}$ and a set of lines $\mathcal{E} \subset E$ to be switched off such that the network congestion of the post-switching network $\Gamma(\GcalE)$ is minimized while satisfying constraints (a), (b), (c) and \eqref{eq:DCPF}:
\begin{align*}
    \min \; & \Gamma(\GcalE),\\
    \text{s.t. } & \text{(a), (b), (c), \eqref{eq:DCPF}}.
\end{align*}
\edit{Unlike TP-PFD, \tppgamma~incorporates the power flow model, making it a more complex optimization problem, as discussed further in \Cref{sec:methods}.}

\subsection{\edit{Alternative tree partitioning problem formulations}}
\label{sec:other-formulations}
The tree partitioning framework presented earlier in this section is rather general and can be tailored further depending on the specific target application. 
Besides the two formulations for minimizing power flow disruption and network congestion presented in the previous two sections, another tree partitioning variant was proposed by \cite{BialekJanusz2021} in which the objective is to minimize the sum of line overloads, that is
$$
r(\calE) = \sum_{(i, j) \in E \setminus \calE} \max\{|f^{(\calE)}_{ij}| - c_{ij}, 0\},
$$
where $f^{(\calE)}_{ij}$ denotes the power flow on line $(i,j)$ on the post-switching network $\GcalE$.

Similarly, one could consider other variants of the tree partitioning problem by considering other objective functions or even by adding additional constraints. 
In particular, we can directly extend any problem formulation from the controlled islanding literature into our tree partitioning framework. 
Recent studies consider, for example, minimizing power imbalance and DC-OPF \cite{Tyuryukanov2022}.

The main methodological contributions in this paper revolve around MILPs, which is not suitable if a nonlinear AC power flow model is considered. 
This limitation also arises for other nonlinear constraints, such as those related to voltage and frequency stability.  
However, as shown in \cite{Atat2022}, one could resort to \edit{Benders'} decomposition approach to address these constraints, but this is outside the scope of this paper.
%%% Local Variables: ***
%%% mode:latex ***
%%% TeX-master: "main.tex"  ***
%%% End: ***
\section{Solution methods}
\label{sec:methods}
In this section, we present two different approaches to solving tree partitioning problems. 
The first approach is to formulate the tree partitioning problem as a MILP, which can be solved optimally. 
The second approach decomposes the tree partitioning problem into two subsequent optimization problems, which results in faster but possibly sub-optimal solutions. 
We refer to the first and second approaches as \textit{single-stage} and \textit{two-stage} approaches, respectively.

\subsection{Single-stage approach}
\label{sec:milp}
The combinatorial nature of the general tree partitioning problem lends itself well to MILP.
In particular, as we only consider linear objective functions and constraints in \tppfd~and~\tppgamma, we can formulate these problems as an MILP and solve them to optimality without the need to resort to heuristic methods. 

We first present a general MILP formulation for tree partitioning problems. 
In the rest of the paper, we refer to any line that is switched off (i.e., it belongs to the subset $\calE$) as \emph{inactive}, and \emph{active} otherwise. 
Let $x_{il} \in \binaryset$ denote whether bus $i\in V$ belongs to cluster $l \in K$. 
Let $y_{ijl} \in \binaryset$ denote whether line $(i,j) \in E$ is an internal edge of cluster $l \in K$, \ie, both endpoints $i$ and $j$ belong to cluster $l$. 
Let $z_{ij} \in \binaryset$ denote whether $(i,j)\in E$ is an \textit{active line} and let $w_{ij} \in \binaryset$ denote whether line $(i,j) \in E$ is an \textit{active cross edge}. 
Finally, let $q_{ij} \in \mathbb{R}$ represent the commodity flow on line $(i,j) \in E$. 

\noindent The general tree partitioning problem can be formulated as%
\begin{subequations}
\label{opt:milp}
\begin{align}
    \min \; & r(\calE),  \label{eq:general-objective} \\
    \text{s.t. } &  x_{il} = 1, &&\quad \forall l \in K, \forall i \in \calG_l, \label{eq:generator_coherency} \\
    &  \sum_{l=1}^{k} x_{il} = 1, &&\quad \forall i \in V,\label{eq:exactly1_partition} \\
    &  y_{ijl} \leq x_{il}, &&\quad \forall (i, j) \in E, \forall l \in K, \label{eq:and1}\\
    &  y_{ijl} \leq x_{jl}, &&\quad \forall (i, j) \in E, \forall l \in K,\label{eq:and2}\\
    &  y_{ijl} \geq x_{il} + x_{jl} - 1, &&\quad \forall (i, j) \in E, \forall l \in K,\label{eq:and3}\\
    &  \sum_{(i, j) \in E} q_{ij} - \sum_{(j, i) \in E} q_{ji} = s_i, &&\quad \forall i \in V, \label{eq:scf2}\\
    &  q_{ij} \geq -(n-1) z_{ij} , &&\quad \forall (i, j) \in E, \label{eq:scf3}\\
    &  q_{ij} \leq (n-1) z_{ij}, &&\quad \forall (i, j) \in E, \label{eq:scf4}\\
    &  \sum_{l=1}^{k} y_{ijl} + w_{ij} = z_{ij}, &&\quad \forall (i,j) \in E, \label{eq:activated_edge}\\
    &  \sum_{(i, j) \in E} w_{ij} = k-1, \label{eq:st} \\
    &  x_{il} \in \binaryset, &&\quad \forall i \in V, \forall l \in K, \label{eq:domain-x}\\
    &  y_{ijl} \in \binaryset, &&\quad \forall (i, j) \in E, \forall l \in K, \label{eq:domain-y}\\
    &  w_{ij} \in \binaryset, &&\quad \forall (i, j) \in E, \label{eq:domain-w}\\
    &  z_{ij} \in \binaryset, &&\quad \forall (i, j) \in E, \label{eq:domain-z}\\
    &  q_{ij} \in \mathbb{R}, &&\quad \forall (i, j) \in E. \label{eq:domain-q}
\end{align}
\end{subequations}

The objective function \Cref{eq:general-objective} is to minimize some risk function $r(\calE)$. 
Constraints~\Cref{eq:generator_coherency} group coherent generators within the same cluster. Constraints~\Cref{eq:exactly1_partition} ensure that each bus belongs to exactly one cluster. Constraints~\cref{eq:and1,eq:and2,eq:and3} state that a line $(i,j)$ is an internal edge of cluster $l$ if and only if $i$ and $j$ belong to the same cluster $l$. Note that if a line $(i, j)$ is not an internal edge, \ie, $y_{ijl} =0$ for all $l \in K$, then this readily implies that $(i,j)$ is a cross edge. Constraints~\cref{eq:scf2,eq:scf3,eq:scf4} represent single commodity flow constraints, which ensure that the post-switching network is connected. Here, we set $s_i=-1$ for all $i\in V \setminus \{1\}$ and $s_1=n-1$.
Constraints~\cref{eq:activated_edge} ensures that if line $(i, j)$ is an internal edge, it must always be an active line. Otherwise, we have that $w_{ij} = z_{ij}$, i.e., a line is active if and only if it is an active cross edge.
Constraint~\cref{eq:st} ensures that the post-switching graph contains exactly $k-1$ active cross edges. Together with the single-commodity flow constraints, this guarantees that the identified clusters must be connected in a tree-like manner\edit{, hence guaranteeing that the active cross edges are cut-edges for the network}. 
Finally, \Cref{eq:domain-x,eq:domain-y,eq:domain-w,eq:domain-z,eq:domain-q} define the variable domains.

We now discuss how to modify \eqref{opt:milp} to formulate \tppfd~and \tppgamma.
The MILP formulation for \tppfd~is 
\begin{subequations}
\label{opt:milp-tp-pfd}
\begin{align}
    \min \; & \sum_{(i, j) \in E} \abs{f_{ij}} (1-z_{ij}), \\
    \text{s.t. } &\text{\cref{eq:generator_coherency,eq:exactly1_partition,eq:and1,eq:and2,eq:and3,eq:scf2,eq:scf3,eq:scf4,eq:st,eq:activated_edge,eq:domain-x,eq:domain-y,eq:domain-w,eq:domain-z,eq:domain-q}.}
\end{align}
\end{subequations}
Note that we only need to replace the objective function, which is to minimize sum of the absolute power flows of the inactive lines. 
The single-stage approach is particularly effective for TP-PFD, as we will show in Section~\ref{sec:experiments}. 

To formulate \tppgamma, we introduce a new variable $\gamma \geq 0$ that represents the network congestion. Furthermore, we introduce the power flow variables $f_{ij} \in \mathbb{R}$ for each line $(i,j) \in E$ and bus angle variables $\theta_i \in \left[\theta^{\min}, \theta^{\max} \right] $ for each bus $i \in V$ to represent the DC power flow equations \eqref{eq:DCPF} as constraints.
The MILP formulation for \tppgamma~is then given by
\begin{subequations}
\label{opt:milp-tp-gamma}
\begin{align}
    \min \; & \gamma,\\
\text{s.t. } &\text{\cref{eq:generator_coherency,eq:exactly1_partition,eq:and1,eq:and2,eq:and3,eq:scf2,eq:scf3,eq:scf4,eq:st,eq:activated_edge,eq:domain-x,eq:domain-y,eq:domain-w,eq:domain-z,eq:domain-q}},\\
    &\gamma \geq \abs{f_{ij}} / c_{ij}, &   \forall \, (i, j) \in E, \label{eq:congestion}\\
    &\sum_{(i,j) \in E} f_{ij} - \sum_{(j, i) \in E} f_{ji} = p_i, & \forall \, i \in V, \label{eq:flow_conservation}\\
    &f_{ij} \leq b_{ij}(\theta_i - \theta_j) + M_{ij}(1-z_{ij}), &  \forall \, (i, j) \in E, \label{eq:dcopf_switch_1}\\
    &f_{ij} \geq b_{ij}(\theta_i - \theta_j) - M_{ij}(1-z_{ij}), &  \forall \, (i, j) \in E, \label{eq:dcopf_switch_2}\\
    &f_{ij} \leq M_{ij}z_{ij}, &   \forall \, (i, j) \in E, \label{eq:dcopf_switch_3}\\
    &f_{ij} \geq -M_{ij}z_{ij}, &   \forall \, (i, j) \in E, \label{eq:dcopf_switch_4}\\
    &\gamma \geq 0, \label{eq:domain-Gamma}\\
    &f_{ij} \in \mathbb{R}, &   \forall \, (i, j) \in E, \label{eq:domain-f}\\
    &\theta_i \in \left[\theta^{\min}, \theta^{\max} \right], &   \forall \, i \in V. \label{eq:domain-theta}
\end{align}
\end{subequations}
Constraints~\cref{eq:congestion} bounds the network congestion from below by each line congestion level. 
Constraints~\cref{eq:flow_conservation} ensures flow conservation at each bus. 
Constraints~\cref{eq:dcopf_switch_1,eq:dcopf_switch_2,eq:dcopf_switch_3,eq:dcopf_switch_4} represent the DC power flow constraints with switching actions. A valid big-M value is given by $M_{ij} = b_{ij}\abs{\theta^{\max} - \theta^{\min}}$ for each line $(i,j)\in E$. 
Finally, \Cref{eq:domain-Gamma,eq:domain-f,eq:domain-theta} represent the variable domains.

\begin{figure*}[!ht]
    \centering
    % Answer: [trim={left bottom right top},clip]
    \subfloat[]{\includegraphics[trim={5cm 2 5cm 2},clip,width=0.30\linewidth]{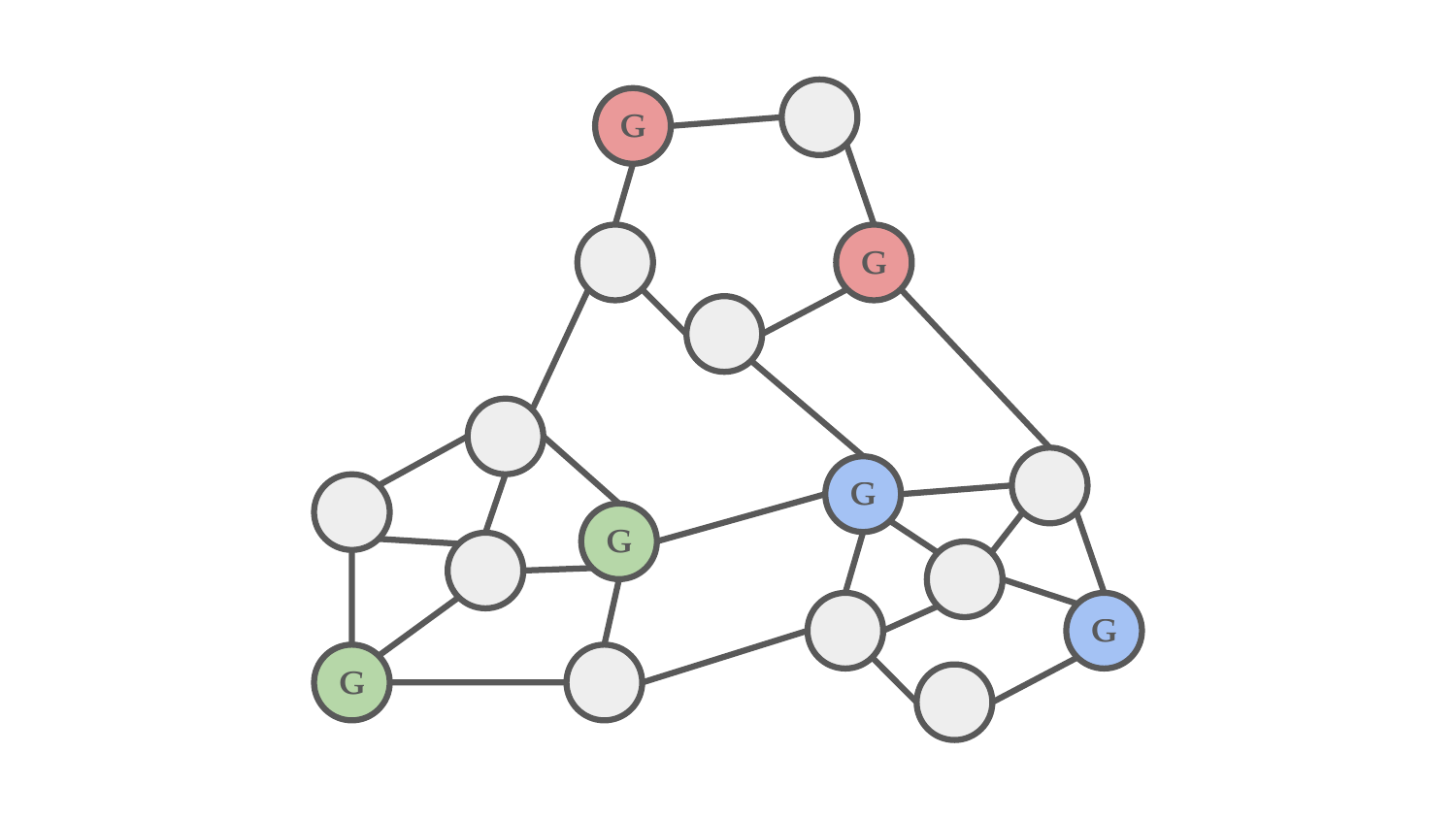}%
    \label{fig_first_case}}
    \hfill
    \subfloat[]{\includegraphics[trim={5cm 2 5cm 2},clip,width=0.30\linewidth]{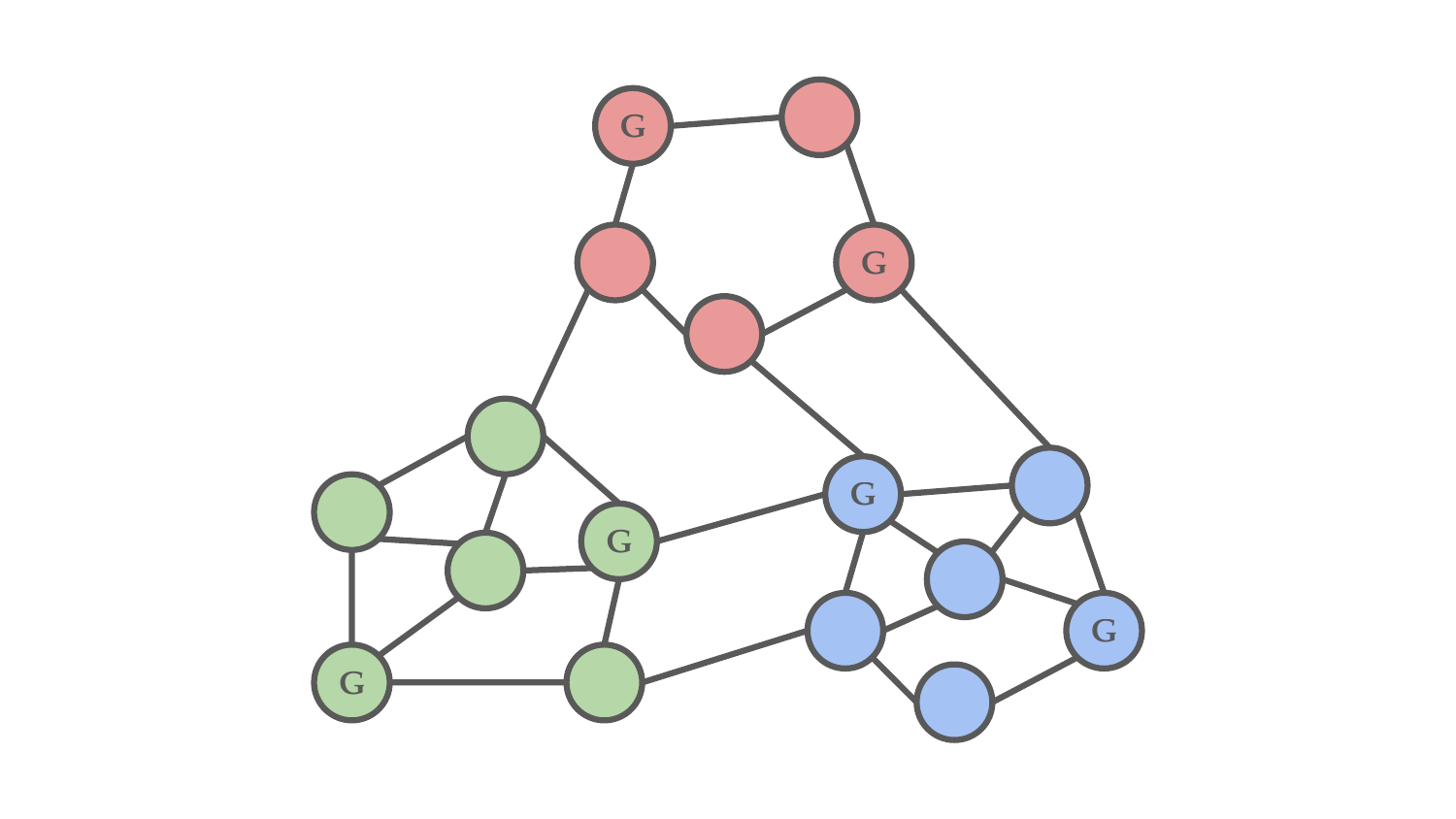}%
    \label{fig_second_case}}
    \hfill
    \subfloat[]{\includegraphics[trim={5cm 2 5cm 2},clip,width=0.30\linewidth]{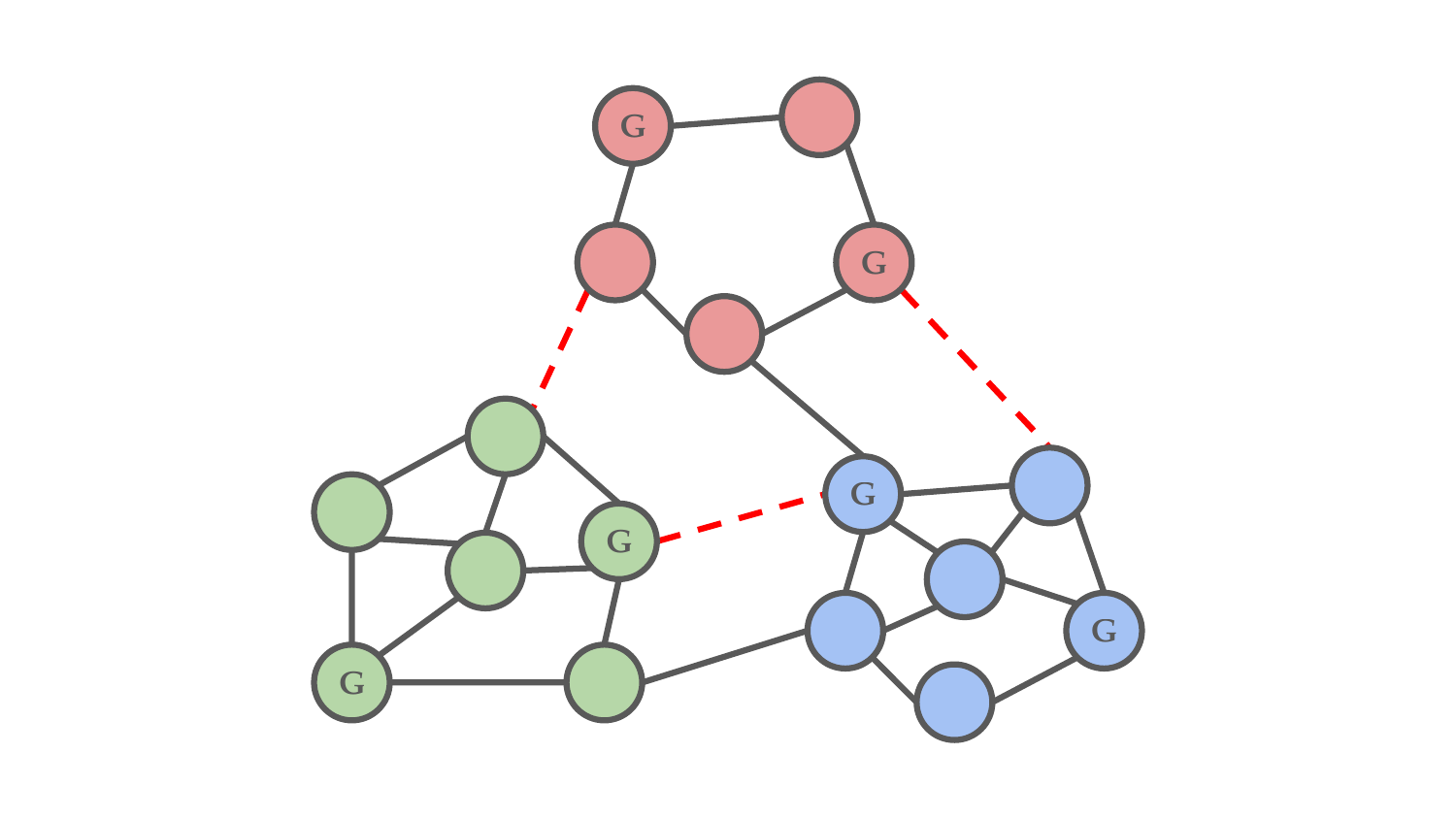}%
    \label{fig_third_case}}
    \caption{
    Illustration of the two-stage approach. (a) A power network $G$ with three generator groups. (b) A 3-partition $\calP$ respecting the generator coherency constraints. (c) A tree partition $\calP$ of $G^\calE$, where the red lines $\calE$ are switched off.}
    \label{fig:two-stage}
    \vspace{-0.45cm}
\end{figure*}

\subsection{Two-stage approach}
\label{sec:two-stage}
Tree partitioning problems can be naturally decomposed into two consecutive stages, which are illustrated in Figure~\ref{fig:two-stage}. In the first stage, one identifies a partition $\calP$ that will serve as the candidate tree partition. In the second stage, one selects the optimal subset of cross edges $\calE \subset E_C(\calP)$ whose removal turns $\calP$ into a tree partition. 

The two-stage approach is useful when the single-stage approach becomes intractable. 
On the other hand, such a heuristic approach may lead to sub-optimal results: there is no guarantee that the identified partition leads to the best switching actions in general.
Nevertheless, we may obtain reasonable results with the advantage of having lower runtimes. We will now discuss the two stages separately in more detail.

\subsubsection{First stage -- Tree Partition Identification (TPI)}
The first stage, which we call the \textit{Tree Partition Identification (TPI)} problem, aims to find a $k$-partition of the power network that respects the generator coherency constraints and will serve as desired tree partition on the post-switching network.

The TPI problem gives rise to the heuristic nature of the two-stage approach: namely, it is hard to say in advance what the optimal partition is with respect to truly optimal line-switching actions. 
However, switching off too many lines or lines with large power flows often leads to large power flow redistribution and results in severe congestion in the resulting network. 
We thus seek a partition such that the resulting line-switching actions have minimal impact on the network. 
One can make several choices for the objective function of TPI, depending on the considered formulation. 
As we are aiming to minimize power flow disruption in TP-PFD, a logical choice is to also minimize power flow disruption for the TPI problem. Furthermore, having low-weight cross edges leads to low power flow redistributions, making this objective function a good candidate for the \tppgamma~problem. 

The MILP for TPI is very similar to \eqref{opt:milp-tp-pfd}, see Appendix~\ref{app:milp-tpi} for its formulation.
The main differences are (i) the exclusion of activated cross edges since all cross edges are deactivated in TPI and (ii) a modification of the single commodity flow constraints to ensure connectivity within each cluster.

One may observe that the TPI problem is identical to the \edit{controlled islanding} problem. 
An alternative objective function that is often considered in the \edit{controlled islanding} literature is the minimization of power imbalance. However, based on our numerical experiments, this objective function is not appropriate for the considered tree partitioning problems since the resulting partition may contain cross edges with very high power flows.

Another widely used method in the \edit{controlled islanding} literature is \textit{spectral clustering}. This is also used in the two-stage tree partitioning method by \cite{BialekJanusz2021,Zocca}. 
However, we remark neither of the aforementioned studies considers the inclusion of generator coherency constraints in their spectral clustering method.

\subsubsection{Second stage -- Optimal Line Switching (OLS)}
 Having identified a good candidate partition $\calP$ in the first stage, the next step is to determine a set of lines $\calE$ to be switched off such that $\calP$ becomes a tree partition of the post-switching network $G^\calE$. We formulate the second stage as an optimization problem, named the \emph{Optimal Line Switching (OLS)} problem. Given a power network $G=(V, E)$ and a $k$-partition $\calP$, the goal of the OLS problem is to remove a subset of cross edges $\calE \subset E_C(\calP)$ to minimize some risk function $r(\calE)$ under the constraint that $\calP$ is a tree partition of $G^\calE$.

An alternative formulation of the OLS problem uses the notion of the reduced graph. Recall that the reduced graph of $G$ given a $k$-partition $\calP$ is defined as $G_{\calP}=(K, E_{C}(\calP))$, i.e., as the graph whose vertices are the clusters of $\calP$ indexed from $1$ to $k$ and whose edges are the cross edges between them. Solving the OLS problem is thus equivalent to computing a spanning tree $T$ on the reduced graph $G_\calP$, such that the removal of lines $\calE = E_C(\calP) \setminus T$ minimizes the risk function $r(\calE)$ on the post-switching network $G^\calE$.

For TP-PFD, the OLS problem is equivalent to the \emph{maximum spanning tree} problem \cite{BialekJanusz2021}. Observe that the spanning tree with maximum weight implies that the set of switched-off lines $\calE$ is of minimum weight. As minimum spanning tree problems can be solved optimally in polynomial time, we can use the negative of the absolute line weights and solve the OLS problem efficiently as well.

For \tppgamma, the OLS problem is more involved as it does not reduce to a polynomial-time solvable spanning tree problem. It is particularly difficult since for any given subset $\calE$, we need to recalculate the power flows on the post-switching graph $G^\calE$ in order to obtain the network congestion $\Gamma(\GcalE)$. In \cite{Zocca}, the OLS problem was solved using a brute-force algorithm that enumerates all possible spanning trees. However, since their number is exponential in the number of clusters $k$, any brute force formulation is intractable for large instances. Instead, we formulate the OLS problem as MILP and solve it as such. We refer to Appendix~\ref{app:milp-ols} for a description of the MILP formulation.

%%% Local Variables: ***
%%% mode:latex ***
%%% TeX-master: "main.tex"  ***
%%% End: ***
\section{Numerical experiments}
\label{sec:experiments}
In this section, we describe our numerical experiments and discuss the results. More specifically, in \Cref{sec:results-pfd} and \Cref{sec:results-gamma}, we compare the performance between the single-stage approach (1-ST) and two-stage approach (2-ST) in solving \tppfd~and \tppgamma, respectively, while in \Cref{sec:cfa}, we run DC cascading failure simulations to compare the effectiveness of different tree partitioning approaches against cascading failures. 

In our numerical experiments, we use a subset of test cases from the \textsc{PGLib-OPF} library \cite{Babaeinejadsarookolaee2021}. From each test case, we extract a power network $G=(V, E)$, where the power injections $\bm{p}$ and flows $\bm{f}$ are computed by solving a DC-OPF problem. 

We obtained generator groups that result in feasible tree partitions using the following procedure. 
We compute a minimum spanning tree for each instance, using the negative absolute power flows as edge weights. 
This spanning tree is iteratively split into $k$ sub-trees. 
At each iteration, we select the largest sub-tree and divide it into two parts such that the ratio between the number of generators in either sub-tree is as close to one as possible.

We ran all experiments on Intel Xeon Gold 6130 CPU processors using 16 cores. 
Our code is implemented in Python, and we use the commercial software Gurobi 9.5.2 to solve MILPs. 
The implementation is openly available at \href{https://github.com/leonlan/tree-partitioning}{https://github.com/leonlan/tree-partitioning}.

For solving an instance with the single-stage approach, we impose a time limit of 600 seconds for the MILP solver.
When the runtime of the corresponding result is 600 seconds, this indicates that the solution may be sub-optimal.
For the two-stage approach, we set a 300-second time limit for solving each of the two separate stages. 

\subsection{Results for \tppfd}
\label{sec:results-pfd}
This section presents the results related to the TP-PFD variant as presented in Section~\ref{sec:transient-stability}.
Table~\ref{tab:results-pfd} presents the computed power flow disruption and runtime of both methods for various power networks and values of $k$. 
We also report the percentage difference in power flow disruption, comparing 2-ST with respect to 1-ST. 
The results show that 1-ST produced optimal solutions for all but the two largest instances. 2-ST obtained optimal solutions in 22 out of 40 instances. 
For the other instances, 14 were solved with a gap of less than 15\%, whereas the remaining 4 instances achieved a gap of larger than 34\%. 
The runtimes of 1-ST were within several seconds in 34 out of 40 instances but increased drastically for larger problem instances.
In contrast, runtimes of 2-ST seem to scale better: the longest measured runtime was only 31 seconds.
In summary, the results show that 1-ST is highly effective at solving TP-PFD optimally for the considered instances.
Moreover, 2-ST obtained optimal results for at least half of the instances and can be considered a fast alternative to 1-ST to produce high-quality solutions for larger instances.

\begin{table}[!ht]
\centering
\setlength{\tabcolsep}{0.57em} % for the horizontal padding
\caption{Results for tree partitioning with as objective minimizing the power flow disruption.}
\label{tab:results-pfd}
\centering
\setlength{\tabcolsep}{0.8em} % for the horizontal padding
\begin{tabular}{llrrrrr}
\toprule
\multicolumn{2}{l}{} & \multicolumn{3}{c}{Objective value (PFD)} & \multicolumn{2}{c}{Runtime (s)} \\
\cmidrule{3-5} \cmidrule{6-7}
 Name        & $k$ \, &    1-ST &  2-ST & \% gap &   1-ST &  2-ST \\
\midrule
EPRI-39 & 2 &                    50 &    50 &    0.00 &    0.09 &    0.06 \\
         & 3 &                    50 &    50 &    0.00 &    0.12 &    0.07 \\
         & 4 &                    50 &    67 &   34.00 &    0.08 &    0.09 \\
         & 5 &                    34 &    67 &   97.06 &    0.13 &    0.12 \\
\midrule
IEEE-57 & 2 &                   158 &   158 &    0.00 &    0.12 &    0.10 \\
         & 3 &                   155 &   155 &    0.00 &    0.17 &    0.12 \\
         & 4 &                   172 &   172 &    0.00 &    0.20 &    0.15 \\
         & 5 &                   172 &   172 &    0.00 &    0.26 &    0.25 \\
\midrule
IEEE-118 & 2 &                   267 &   267 &    0.00 &    0.22 &    0.25 \\
         & 3 &                   277 &   277 &    0.00 &    0.23 &    0.24 \\
         & 4 &                   786 &   786 &    0.00 &    0.27 &    0.29 \\
         & 5 &                   812 &   812 &    0.00 &    0.31 &    0.42 \\
\midrule
GOC-179 & 2 &                   252 &   252 &    0.00 &    0.29 &    0.27 \\
         & 3 &                  1944 &  1944 &    0.00 &    1.04 &    0.44 \\
         & 4 &                  2796 &  2796 &    0.00 &    5.40 &    0.55 \\
         & 5 &                  2796 &  2796 &    0.00 &   66.74 &    0.72 \\
\midrule
IEEE-300 & 2 &                   193 &   193 &    0.00 &    0.45 &    0.51 \\
         & 3 &                   312 &   312 &    0.00 &    0.54 &    0.66 \\
         & 4 &                   909 &   909 &    0.00 &    0.63 &    0.80 \\
         & 5 &                  1006 &  1148 &   14.12 &    0.78 &    1.05 \\
\midrule
GOC-500 & 2 &                   560 &   560 &    0.00 &    0.74 &    0.99 \\
         & 3 &                   740 &   800 &    8.11 &    0.96 &    1.32 \\
         & 4 &                  1221 &  1281 &    4.91 &    1.65 &    1.60 \\
         & 5 &                  1236 &  1381 &   11.73 &    2.61 &    1.92 \\
\midrule
SDET-588 & 2 &                   135 &   135 &    0.00 &    0.73 &    0.96 \\
         & 3 &                   436 &   436 &    0.00 &    0.80 &    1.24 \\
         & 4 &                   561 &   561 &    0.00 &    1.03 &    1.56 \\
         & 5 &                   568 &   768 &   35.21 &    1.22 &    1.77 \\
\midrule
GOC-793 & 2 &                   673 &   673 &    0.00 &    0.93 &    1.28 \\
         & 3 &                   917 &   975 &    6.32 &    1.27 &    1.63 \\
         & 4 &                   917 &  1030 &   12.32 &    2.19 &    4.44 \\
         & 5 &                  1048 &  1480 &   41.22 &    3.31 &   15.04 \\
\midrule
RTE-1888 & 2 &                   788 &   835 &    5.96 &    2.80 &    4.27 \\
         & 3 &                  1623 &  1670 &    2.90 &    3.80 &    5.24 \\
         & 4 &                  3757 &  3804 &    1.25 &   48.67 &    6.18 \\
         & 5 &                  5245 &  5361 &    2.21 &  600.00 &   30.84 \\
\midrule
RTE-2848 & 2 &                   889 &   955 &    7.42 &    4.26 &    6.53 \\
         & 3 &                  1624 &  1690 &    4.06 &  167.35 &    8.13 \\
         & 4 &                  2259 &  2286 &    1.20 &  388.56 &    9.77 \\
         & 5 &                  3197 &  3224 &    0.84 &  600.00 &   13.25 \\
\bottomrule
\end{tabular}

\end{table}

\subsection{Results for \tppgamma}
\label{sec:results-gamma}
This section presents the results related to the \tppgamma~variant as presented in Section~\ref{sec:network-congestion}.
We remark that \tppgamma~is significantly harder to solve than \tppfd: especially for large power networks and larger values of $k$, \tppgamma~could not be solved with the single-stage approach using the 600-second time limit.
We, therefore, decided to use a \textit{warm-started} solution for all instances larger than GOC-179.
In particular, we use the solution obtained from solving \tppfd~using 1-ST and use this as the initial solution when solving \tppgamma~with 1-ST. 
This approach makes a direct comparison between the 1-ST and 2-ST approaches unfair, but it allows us to find (near-)optimal solutions to evaluate the solution quality of 2-ST. 

Table~\ref{tab:results-mc} reports the results from running 1-ST and 2-ST on the considered power network instances. 
All networks start with network congestion of 1, which is a direct result of initializing the networks with DC-OPF.

We first describe the results related to the objective value, i.e., the network congestion. 
The warm-started version of 1-ST obtained optimal solutions for 30 out of 40 instances.
The 2-ST approach produced optimal solutions for 20 out of the 40 instances. 
For the remaining instances, 10 were solved with a gap of at most 16\%, and the other 10 were solved with a gap larger than 26\%.
It was possible to find an optimal solution resulting in no congested lines for 26 out of 40 instances.
The remaining 14 instances could not be solved without creating congested lines. However, in 8 out of these 14 instances, the network congestion remained quite reasonable, below 1.10. Only for 6 instances the network congestion exceeds 1.58. 

In terms of runtimes, we only discuss the results of 2-ST. Overall, the runtimes of 2-ST are similar to those for solving the TP-PFD problem with 2-ST.
This is due to the TPI problem being identical, and most of the time is spent on solving the TPI problem. 
The difference in the runtimes between TP-PFD and \tppgamma~is due to solving the OLS problem. 
This problem is solved very fast: it took at most 5 seconds, with the exception of RTE-1888 and $k=4$, for which it took about 18 seconds. 

In summary, the results for \tppgamma~show that it is possible to partition a network optimally with low congestion in most cases. 
However, minimizing the network congestion is much harder than minimizing the power flow disruption, and consequently, the single-stage approach is not a viable method of solving \tppgamma.
Similar to the results for TP-PFD, the two-stage method runs faster, but the resulting solutions are often of lesser quality.
One main difficulty when minimizing network congestion is that line congestion is a local property: even if many of the switched lines belong to the optimal set, there may be a single line switching action that can cause extremely high congestion on a nearby line. 
As a result, a partition with minimal power flow disruption may not be optimal with respect to the subsequent line-switching actions, since it does not explicitly take into account this congestion.

\begin{table}[!ht]
\centering
\setlength{\tabcolsep}{0.57em} % for the horizontal padding
\caption{
Results for tree partitioning with as objective minimizing the network congestion.
}
\label{tab:results-mc}
{\centering
\setlength{\tabcolsep}{0.8em} % for the horizontal padding
\begin{tabular}{lllrrrr}
\toprule
\multicolumn{2}{l}{} & \multicolumn{3}{c}{Objective value ($\Gamma$)} & \multicolumn{2}{c}{Runtime (s)} \\
\cmidrule{3-5} \cmidrule{6-7}
 Name        & $k$ \, &    1-ST &  2-ST & \% gap &   1-ST &  2-ST \\
\midrule
EPRI-39 & 2 &               1.00 &  1.00 &    0.00 &    0.16 &   0.09 \\
         & 3 &               1.00 &  1.00 &    0.00 &    0.49 &   0.11 \\
         & 4 &               1.00 &  1.00 &    0.00 &    0.64 &   0.14 \\
         & 5 &               1.00 &  1.00 &    0.00 &    0.37 &   0.16 \\
\midrule
IEEE-57 & 2 &               0.88 &  1.01 &   14.77 &    1.03 &   0.15 \\
         & 3 &               0.88 &  1.01 &   14.77 &    2.23 &   0.19 \\
         & 4 &               0.88 &  1.02 &   15.91 &    8.45 &   0.32 \\
         & 5 &               0.88 &  1.02 &   15.91 &   11.04 &   0.38 \\
\midrule
IEEE-118 & 2 &               1.00 &  1.00 &    0.00 &    1.64 &   0.32 \\
         & 3 &               1.00 &  1.00 &    0.00 &    4.87 &   0.34 \\
         & 4 &               1.48 &  1.48 &    0.00 &    9.75 &   0.56 \\
         & 5 &               1.48 &  1.48 &    0.00 &  150.38 &   0.72 \\
\midrule
GOC-179 & 2 &               1.04 &  1.07 &    2.88 &   15.40 &   0.44 \\
         & 3 &               1.58 &  1.63 &    3.16 &  268.10 &   0.57 \\
         & 4 &               1.58 &  1.63 &    3.16 &  600.00 &   0.78 \\
         & 5 &               1.42 &  1.63 &   14.79 &  600.00 &   1.17 \\
\midrule
IEEE-300 & 2 &               1.06 &  1.16 &   9.43 &  600.00 &   1.01 \\
         & 3 &               1.09 &  1.16 &   6.42 &  600.00 &   1.26 \\
         & 4 &               1.09 &  1.38 &  26.61 &  600.00 &   1.64 \\
         & 5 &               1.09 &  1.63 &  49.54 &  600.00 &   2.16 \\
\midrule
GOC-500 & 2 &               1.00 &  1.00 &   0.00 &    2.11 &   1.35 \\
         & 3 &               1.00 &  1.00 &   0.00 &    2.52 &   1.74 \\
         & 4 &               1.00 &  1.00 &   0.00 &    3.59 &   2.34 \\
         & 5 &               1.00 &  1.00 &   0.00 &    4.91 &   2.84 \\
\midrule
SDET-588 & 2 &               1.02 &  1.44 &  41.18 &  204.12 &   1.53 \\
         & 3 &               1.00 &  1.32 &  32.00 &  256.37 &   2.12 \\
         & 4 &               1.00 &  1.29 &  29.00 &  182.36 &   2.56 \\
         & 5 &               1.00 &  1.26 &  26.00 &   19.10 &   3.15 \\
\midrule
GOC-793 & 2 &               1.05 &  2.06 &  96.19 &  600.00 &   2.22 \\
         & 3 &               1.00 &  1.57 &  57.00 &  600.00 &   3.08 \\
         & 4 &               1.01 &  1.57 &  55.45 &  600.00 &   4.93 \\
         & 5 &               1.07 &  1.64 &  53.27 &  600.00 &  14.77 \\
\midrule
RTE-1888 & 2 &               1.00 &  1.00 &   0.00 &    9.56 &   5.49 \\
         & 3 &               1.00 &  1.00 &   0.00 &   11.18 &   6.87 \\
         & 4 &               1.00 &  1.00 &   0.00 &   12.77 &  23.40 \\
         & 5 &               1.00 &  1.00 &   0.00 &   21.15 &  34.96 \\
\midrule
RTE-2848 & 2 &               1.00 &  1.00 &   0.00 &   15.84 &   8.77 \\
         & 3 &               1.00 &  1.00 &   0.00 &   31.00 &  11.43 \\
         & 4 &               1.00 &  1.00 &   0.00 &   48.46 &  13.36 \\
         & 5 &               1.00 &  1.00 &   0.00 &  144.15 &  17.89 \\
\bottomrule
\end{tabular}
}
\end{table}

\subsection{Cascading failure simulations}
\label{sec:cfa}
In this final numerical experiment, we run cascading failure simulations using the DC power flow model to compare the effectiveness of different tree partitioning strategies. 
For a given power network, a single cascading failure simulation is initiated as follows. 
The simulation starts by removing a selected line, which may result in multiple disconnected components.
Within each component, we readjust the power imbalance if necessary, using proportional load shedding or generation curtailing.
More specifically, when a component has more load than generation, then we lower all loads by some proportion until it matches the generation. Otherwise, we curtail the generation with some proportion to match the load.
After adjusting generation and load, we re-run the DC power flow equations and identify the overloaded lines.
If no lines are overloaded or all loads are shed, we stop the cascading failure simulation and register the total lost load.
Otherwise, the overloaded lines are removed, and the cascading failure simulation continues.

For a given power network, we run a number of cascading failure simulations by removing each line once, and we collect the average lost load over all simulations.
As a baseline, we run a cascading failure simulation on the original power network.
We compare this with the power networks obtained by solving TP-PFD and \tppgamma~using the original power network \edit{and applying 
both the single-stage approach and two-stage approach for $k=2,3,4$.}
After tree partitioning, we run DC-OPF on the network and then run the cascading failure simulations.
Results are labeled as PFD (1-ST), PFD (2-ST), NC (1-ST), and NC (2-ST).

We only consider a subset of the power networks that meet the following criteria. 
First, the average lost load of the original network during the cascading failure simulations must be more than 1\% of the total network load. 
This excludes the networks RTE-1888 and RTE-2848. 
Second, the networks after tree partitioning must have complete results for the cascading failure simulations. 
This criterion excludes IEEE-118 because running DC-OPF for $k=4$ is infeasible after tree partitioning. 

\begin{figure}
    \centering
    \includegraphics[width=0.5\textwidth]{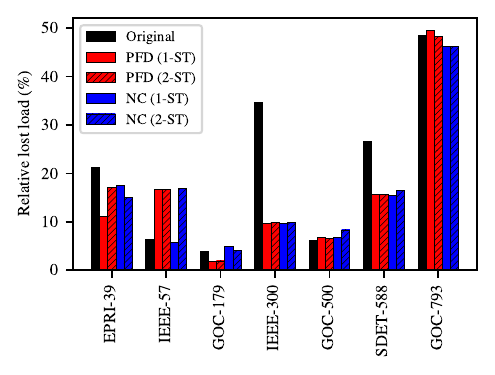}
    \caption{
    Average lost loss (as \% of initial network load) during cascading failure simulations. Results compare the original network against four tree-partitioned networks: power flow disruption (PFD) and network congestion (NC), each solved with single-stage (1-ST) and two-stage (2-ST) approaches.
    }
    \label{fig:results-cfs}
\end{figure}

\Cref{fig:results-cfs} shows the results of the cascading failure simulations, where the values represent the average lost load relative to the total load of the original network before cascading failures.
The results for \edit{tree-partitioned networks} are averaged over all considered values of $k=2,3,4$.
\edit{
Our simulations reveal that, on average, tree partitioning reduces the lost load due to cascading failures when compared to the original network.
The original network shows 21.1\% average lost load, while the tree-partitioned networks perform better: PFD (1-ST) at 15.9\%, PFD (2-ST) at 16.5\%, NC (1-ST) at 15.3\%, and NC (2-ST) at 16.7\%. 
The single-stage approach slightly outperforms two-stage methods, and NC (1-ST) is marginally better than PFD (1-ST), though the differences are small.}
When looking specifically at the individual test cases, tree partitioning shows a large reduction in lost load for EPRI-39, IEEE-300, and SDET-588. In particular, the lost load is reduced from $35\%$ to $10\%$ for the IEEE-300 network.
Mixed results appear for IEEE-57, GOC-179, and GOC-793: \edit{here, some tree partitioning methods perform slightly better than the original network, whereas others perform slightly worse.}
Only GOC-500 shows no benefit from tree partitioning: \edit{in all cases, the lost load is slightly higher than for the original network, with the maximum difference being 2.2\% for NC (2-ST).}

\edit{The results demonstrate that tree partitioning can help to reduce lost load in cascading failures.
However, further analysis and investigation are required to thoroughly understand the impact, even when tree partitioning cannot result in a stable network, like in the case of IEEE-118.}

%%% Local Variables: ***
%%% mode:latex ***
%%% TeX-master: "main.tex"  ***
%%% End: ***
\section{Conclusion}
\label{sec:conclusion}
This work explores tree partitioning strategies in transmission networks, which are emerging as a valid alternative emergency measure to mitigate cascading failures. 
The key idea is altering the network topology strategically via line-switching actions to create a more hierarchical structure that prevents the global propagation of failures.

In this paper, we present a new comprehensive and unified framework based on MILP for tree partitioning problems and show how it can encompass many variants already proposed in the literature.
In particular, it can accommodate various objective functions, among which power flow disruption and network congestion, and include constraints related to grouping coherent generators.
Moreover, we provide extensive numerical results for various tree partitioning problems, comparing results between our proposed exact formulation and two-stage heuristic.
We demonstrate that the exact approach achieves better solutions while the heuristic approach obtains good solutions with considerably lower runtimes.
Lastly, by means of extensive cascading failure simulations, we compare the reduction in lost load with respect to the original topologies that various tree partitioning strategies can achieve.

Although our two-stage heuristic approach has shown promise in enhancing scalability, we believe its performance can be enhanced even further by improving the quality of the network partition identified in the first stage, either by including more engineering details or by using more advanced ideas proposed in the study of complex networks.

Furthermore, exploring tree partitioning strategies using the AC power flow model was beyond the scope of this paper, whose core focus was providing a unified optimization framework using solely linear constraints. 
In future work, we hope to account for this more realistic power flow model, as well as \edit{including generation and load adjustments as decision variables}.
\appendices
\section{MILP formulation for TPI}
\label{app:milp-tpi}
The MILP formulation for TPI follows largely the same constraints as \cref{opt:milp-tp-pfd}. 
The main differences are that cross edges are always inactive, and the single commodity flow constraints must be slightly adjusted to accommodate connectivity within clusters.
For each generator group $\mathcal{G}_l, l \in K$, we select one generator bus that represents the source node for the corresponding cluster.
Define $s_i=n-1$ if $i \in V$ is a source bus and $s_{i} = -1$ otherwise.
The MILP formulation of TPI is then given by
\begin{subequations}
\begin{align}
\label{opt:tpi}
    \min \; &\sum_{(i, j) \in E} \abs{f_{ij}} (1-z_{ij}), \\
    \text{s.t. } & \text{\cref{eq:generator_coherency,eq:exactly1_partition,eq:and1,eq:and2,eq:and3,eq:scf3,eq:scf4}},\\
    & \text{\cref{eq:domain-x,eq:domain-y,eq:domain-z,eq:domain-q}}, \\
    & \sum_{l \in K} y_{ijl} = z_{ij}, &\quad \forall (i, j) \in E,\label{eq:tpi-active-line}\\
    & \sum_{(i, j) \in E} q_{ij} - \sum_{(j, i) \in E} q_{ji} \le s_i, &\quad \forall i \in V. \label{eq:tpi-scf2}
\end{align}
\end{subequations}
Constraints~\cref{eq:tpi-active-line} ensure that internal edges are active while cross edges are inactive. 
Constraints~\cref{eq:tpi-scf2} ensure that source nodes have an outflow of at most $n-1$ commodity units and other nodes have an inflow of at least 1 unit.

\section{MILP formulation for OLS}
\label{app:milp-ols}
The MILP formulation for OLS (as part of solving \tppgamma) uses the same constraints as \cref{opt:milp-tp-gamma}.
However, now that a candidate partition $\calP = \{\calV_1, \dots, \calV_k\}$ has been identified, we can fix all partitioning variables in the MILP. In particular, the MILP for OLS is given by \eqref{opt:milp-tp-gamma} in combination with the following constraint:
\begin{align}
    x_{il} = 1, & \quad \forall l \in K, \forall i \in \calV_l.
\end{align}

\bibliographystyle{IEEEtran}
\bibliography{bibliography}

\begin{IEEEbiography}
[{\includegraphics[width=1in,height=1.25in,clip,keepaspectratio]{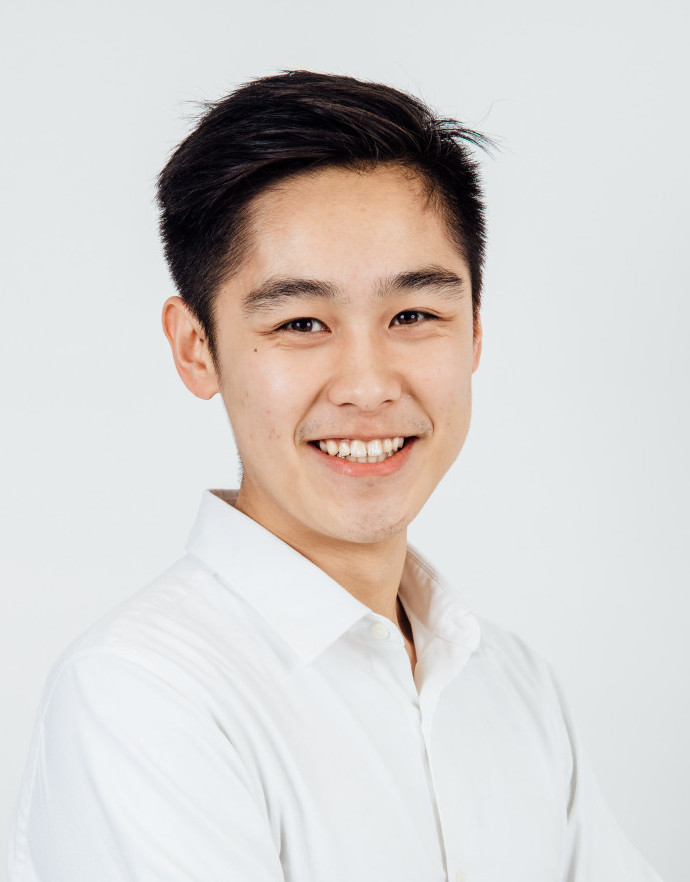}}]
{Leon Lan} received a B.Sc. degree in liberal arts and sciences from the Amsterdam University College, The Netherlands, in 2015 and an M.Sc. degree in operations research from the Vrije Universiteit Amsterdam, The Netherlands.
Since 2021, he has been a Ph.D. student at the Department of Mathematics, Vrije Universiteit Amsterdam, The Netherlands.
His research focuses on the integration of production scheduling and vehicle routing in large-scale supply chains.

\end{IEEEbiography}

\begin{IEEEbiography}
[{\includegraphics[width=1in,height=1.25in,clip,keepaspectratio]{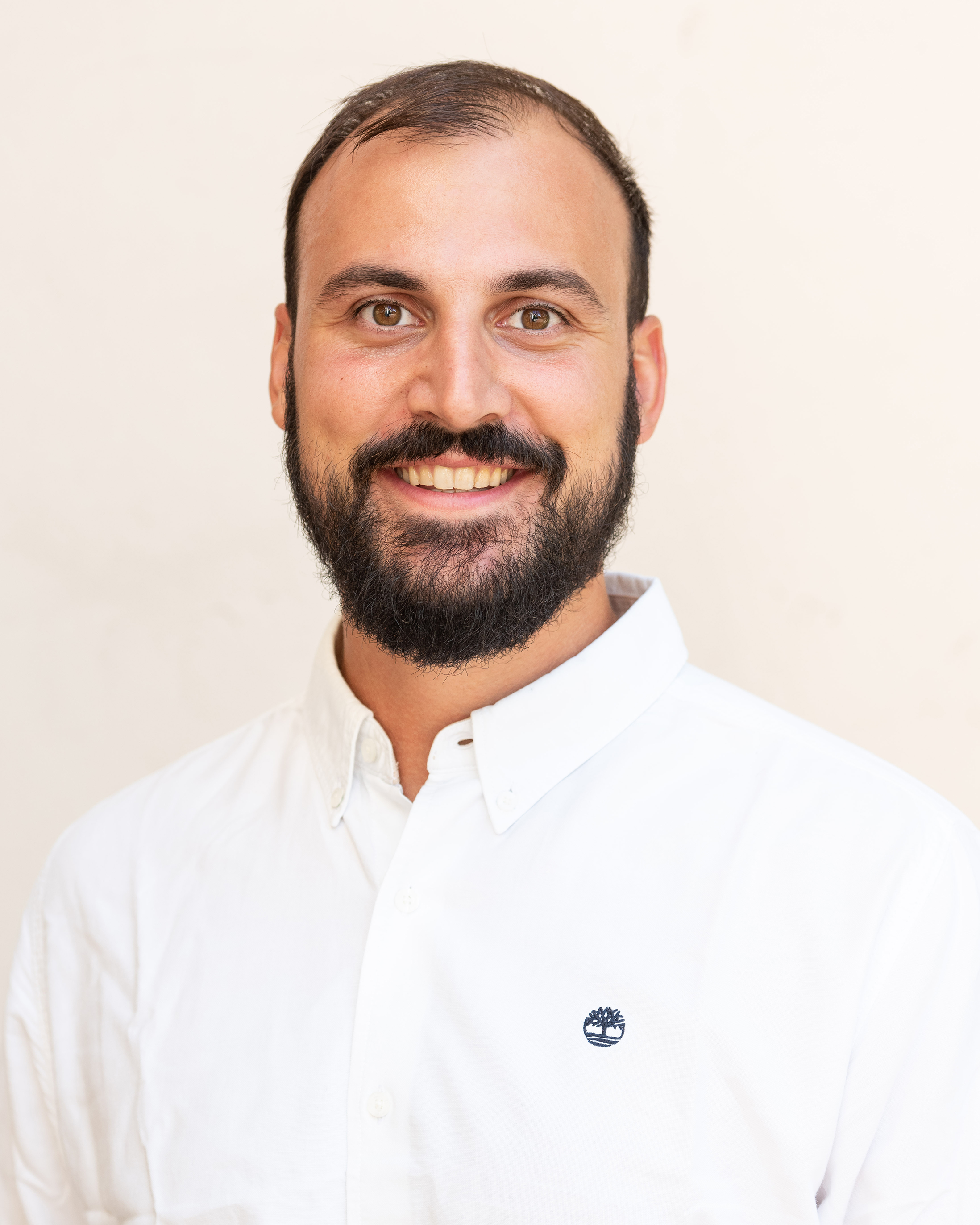}}]
{Alessandro Zocca} received his Ph.D. degree in mathematics from the University of Eindhoven, in 2015. He has been a postdoctoral researcher at CWI Amsterdam, and then at the California Institute of Technology, where he was supported by his personal NWO Rubicon grant. Since 2019, he is a tenure-track Assistant Professor with the Department of Mathematics of the Vrije Universiteit Amsterdam. His work lies mostly in the area of applied probability, learning, and optimization, drawing motivation from applications to power systems reliability.
\end{IEEEbiography}

%\newpage



%\vfill


% \input{reviews}

\end{document}